\documentclass[12pt]{article}

\usepackage{amsmath,amsthm,amscd,amssymb}

\usepackage[mathscr]{eucal}

\usepackage{amssymb}

\usepackage{latexsym}

\usepackage{enumerate} 

\newtheorem{Thm}{Theorem}[section]

\newtheorem{Prop}[Thm]{Proposition}

\newtheorem{Cor}[Thm]{Corollary}

\newtheorem{Lem}[Thm]{Lemma}

\theoremstyle{definition}

\newtheorem{Rem}[Thm]{Remark}

\DeclareMathOperator{\re}{Re}

\DeclareMathOperator{\tr}{tr}

\DeclareMathOperator{\rank}{rank}

\DeclareMathOperator{\ord}{ord}

\DeclareMathOperator{\diag}{diag} 

\numberwithin{equation}{section}

\newcommand{\Q}{\mathbb{Q}}

\newcommand{\C}{\mathbb{C}}

\newcommand{\Z}{\mathbb{Z}}

\newcommand{\F}{\mathbb{F}}

\newcommand{\gamt}{\Gamma_{2}}

\newcommand{\trn}[1][1]{{}^t \hspace{-#1pt}}

\newcommand{\mat}[4]{\begin{pmatrix} #1 & #2 \\ #3 & #4 \end{pmatrix}}

\newcommand{\smat}[4]{\begin{smallmatrix} #1 & #2 \\ #3 & #4 \end{smallmatrix}}

\begin{document}

\renewcommand{\S}{\ensuremath{\mathchar"0278 }}

\title{A ring of symmetric Hermitian modular forms of degree $2$ with integral Fourier coefficients~II}
\author{Toshiyuki Kikuta}
\maketitle


\begin{abstract}
We exhibit 24 homogeneous generators for the graded \(\Z\)-algebra of
symmetric Hermitian modular forms of degree 2 over \(\Q(i)\) with integral
Fourier coefficients, where the weight \(k\) is even and the character is
\(\det^{k/2}\). We prove that 24 is the minimum number of algebra generators.
For the subalgebra in weights divisible by 4, Kurihara reduced the
generating set given in our earlier work to 11 forms. We show that 13 further
forms generate the direct sum of the components in weights congruent
to 2 modulo 4 as a module over this subalgebra.
For \(p=2\) and \(p=3\), we show that the ring of Fourier expansions modulo
\(p\) is a polynomial ring over \(\F_p\) generated by the reductions of
4 forms of weights \(8,12,12,16\).
We also record the diagonal Sturm bound \(\lfloor k/8\rfloor\) modulo
every prime \(p\) and verify its sharpness in every even weight \(k\ne2\).
\end{abstract}

\section{Introduction}\label{sec:ring}

Let \(\mathcal O=\Z[i]\). For each even integer \(k\ge0\), let
\(M_k=M_k(\gamt(\mathcal O);\Z)\) denote the \(\Z\)-module of symmetric
Hermitian modular forms of degree 2, weight \(k\), and character \(\det^{k/2}\)
on \(\gamt(\mathcal O)\) with integral Fourier coefficients. In odd weight
we have \(M_k=0\), as explained in Section~\ref{sec:herm-fourier}. Thus
the full graded \(\Z\)-algebra is
\[
 B=\bigoplus_{k\ge0}M_k=\bigoplus_{m\ge0}M_{2m}.
\]
We determine a minimal set of homogeneous generators for \(B\) over \(\Z\).

For the subalgebra of Hermitian modular forms in weights divisible by 4, write
\[
 A_{\mathrm H}=\bigoplus_{m\ge0} M_{4m}.
\]
The $24$ forms of \cite[Theorem~3.7]{Kikuta} generate \(A_{\mathrm H}\) over \(\Z\).
Kurihara~\cite{Kurihara} reduced this to $11$ generators.
We extend this result to all even weights and prove that the minimum
number of \(\Z\)-algebra generators of \(B\) is $24$.

\noindent\textit{Notation for the generators.}
The subscript of each generator gives its weight. For \(4\mid k\) we
use \(H_k\), with \(I_k\) for an additional generator in the same weight.
For \(k\equiv 2\pmod{4}\) the corresponding symbols are \(K_k\) and \(L_k\).
Thus weight 12 has \(H_{12},I_{12}\), and weight 30 has \(K_{30},L_{30}\).
We write \(J_k\) for auxiliary forms and \(J_{k,j}\) when an additional index is needed.
The subscripts \(0,2\) on \(\mathcal G_0,\mathcal G_2\) indicate
the weights of their members modulo 4.

\begin{Thm}\label{thm:main}
The graded \(\Z\)-algebra \(B\) is generated by the union of the
following 2 sets.
\[
\begin{aligned}
 &\mathcal{G}_0
  =\{H_4,H_8,H_{12},I_{12},H_{16},H_{20},
       H_{24},H_{28},H_{36},H_{40},H_{48}\},\\
 &\mathcal{G}_2
  =\{K_6,K_{10},K_{14},K_{18},K_{22},K_{26},K_{30},L_{30},K_{34},K_{38},K_{42},K_{46},K_{50}\}.
\end{aligned}
\]
More precisely,
\[
 A_{\mathrm H}=\Z[\mathcal G_0],\qquad
 \bigoplus_{k\equiv2\pmod{4}}M_k
     =\sum_{g\in\mathcal G_2} A_{\mathrm H}\,g.
\]
Every \(\Z\)-algebra generating set of \(B\) has at least 24 members.
\end{Thm}

Section~\ref{sec:prelim} fixes notation for Hermitian and Siegel modular
forms. In Section~\ref{sec:forms} we give the defining formulas and use
known congruences to verify integrality. We also compute the restrictions of the generators in low weights. In Section~\ref{sec:generation}
we construct integral Hermitian lifts and prove that the 24 forms generate
the algebra. Section~\ref{sec:minimality} proves minimality and concludes
with the answer to the question on integral lifts in weight 46. Section~\ref{sec:sturm-and-reduction}
also records the diagonal Sturm bound and its sharpness for every prime. For \(p=2,3\)
we then describe the ring of Fourier expansions modulo \(p\) and compute
the minimum number of generators of \(B\otimes_{\Z}\F_p\).

\section{Preliminaries}\label{sec:prelim}

\subsection{Symmetric Hermitian modular forms and their Fourier expansions}\label{sec:herm-fourier}

A Hermitian modular form of degree 2 is a holomorphic function on the
Hermitian upper half-space
\[
 \mathbb H_2=\{Z\in M_2(\C)\mid (Z-\trn\overline Z)/(2i)>0\}.
\]
For \(K=\Q(i)\) with ring of integers \(\mathcal O=\Z[i]\), the Hermitian
modular group is
\[
 \gamt(\mathcal O)=
 \{M\in\mathrm{GL}_4(\mathcal O)\mid \trn\overline M J_2 M=J_2\},
 \qquad J_2=\left(\smat{0}{-I_2}{I_2}{0}\right).
\]
An element \(M=\left(\smat{A}{B}{C}{D}\right)\) of \(\gamt(\mathcal O)\)
acts on \(\mathbb H_2\) by \(Z\mapsto M\langle Z\rangle=(AZ+B)(CZ+D)^{-1}\).

For each even integer \(k\ge0\), let \(M_k(\gamt(\mathcal O))\) denote
the complex vector space of symmetric Hermitian modular forms of weight
\(k\) and character \(\det^{k/2}\).
Thus a holomorphic function \(F\) on \(\mathbb H_2\) belongs to this
space precisely when \(F(\trn Z)=F(Z)\) and
\[
 F(M\langle Z\rangle)=\det(M)^{k/2}\det(CZ+D)^k F(Z)
 \qquad (M\in\gamt(\mathcal O)).
\]
The scalar \(iI_4\) fixes \(Z\) and has determinant 1. Its lower right
block \(iI_2\) has determinant \(-1\), so \(F(Z)=(-1)^kF(Z)\) for
any form whose character factors through the determinant. This gives
\(M_k=0\) in odd weight.

Every \(F\in M_k(\gamt(\mathcal O))\) is invariant under the translations
\(Z\mapsto Z+B\) by Hermitian matrices \(B\in M_2(\mathcal O)\) and has
a Fourier expansion of the form
\[
 F(Z)=\sum_{0\le H\in\Lambda_2(K)}a_F(H)e^{2\pi i\tr(HZ)},
\]
where \(\Lambda_2(K)\) is the lattice of Hermitian matrices given by
\[
 \Lambda_2(K)=\left\{\mat{m}{(r+si)/2}{(r-si)/2}{n}\ \middle|\
                m,r,s,n\in\Z\right\}.
\]
We write \(H=(m,r,s,n)\) for the displayed matrix and \(a_F(m,r,s,n)\)
for its Fourier coefficient \(a_F(H)\). The condition \(0\le H\) is equivalent
to \(m,n\ge0\) and \(r^2+s^2\le4mn\).
For a subring \(R\subseteq\C\), define
\[
 M_k(\gamt(\mathcal O);R)=
 \bigl\{F\in M_k(\gamt(\mathcal O))\mid
   a_F(H)\in R\quad(0\le\forall H\in\Lambda_2(K))\bigr\}.
\]
As in the introduction, we abbreviate \(M_k(\gamt(\mathcal O);\Z)\)
to \(M_k\).

For \(Z=(z_{ij})\), put \(e(t)=\exp(2\pi i t)\) and
\[
 \begin{aligned}
 \dot q_{11}&=e(z_{11}),&
 \dot q_{22}&=e(z_{22}),\\
 \dot q_{12}&=e\bigl((z_{12}+z_{21})/2\bigr),&
 \ddot q_{12}&=e\bigl((z_{12}-z_{21})/(2i)\bigr).
 \end{aligned}
\]
For \(H=(m,r,s,n)\) we have
\[
 \tr(HZ)=mz_{11}+r(z_{12}+z_{21})/2+s(z_{12}-z_{21})/(2i)+nz_{22}.
\]
Thus the Fourier expansion is
\[
 F(Z)=\sum_{\substack{m,n\in\Z_{\ge0},\ r,s\in\Z\\ r^2+s^2\le4mn}}
 a_F(m,r,s,n)\,
 \dot q_{11}^{m}\dot q_{12}^{r}\ddot q_{12}^{s}\dot q_{22}^{n}.
\]

The definitions of the 2 off-diagonal variables in
\cite[\S2.2]{Kikuta} contain a misprint. With those definitions the
exponents in \(e(\tr(HZ))\) are \((-s,r)\). This gives the same
Fourier coefficients and hence the same results. Indeed, for
\(V=\diag(i,1)\) the transformation law applied to \(\diag(V,V)\)
gives \(F(VZV^{-1})=F(Z)\). Comparing coefficients gives
\[
 a_F(m,r,s,n)=a_F(m,-s,r,n),
\]
since the upper right entry of \(VHV^{-1}\) is \((-s+ri)/2\).

\subsection{Hermitian Eisenstein series and the 5 basic forms}\label{sec:basic-five}

To construct the generators we start with 5 forms obtained from Hermitian Eisenstein series.
Let \(P=\{\left(\begin{smallmatrix}A&B\\0&D\end{smallmatrix}\right)\in\gamt(\mathcal O)\}\).
For every even integer \(k\ge6\), define the Hermitian Eisenstein series of
degree 2 by the absolutely convergent sum
\[
 E_k(Z):=\sum_{M\in P\backslash\gamt(\mathcal O)}
             \det(M)^{k/2}\det(CZ+D)^{-k},
 \qquad M=\left(\begin{smallmatrix}A&B\\ C&D\end{smallmatrix}\right).
\]
In weight 4 we take \(E_4\) given by Krieg's Maa{\ss} lift
\cite[\S8, Remark~(b) after Corollary~B]{Krieg}. These series have constant term 1.
The Fourier coefficients of \(E_k\) can be computed by the explicit
formula of Krieg and Dern~\cite{Dern}.
We define the 5 basic forms \(H_4,K_6,H_8,K_{10},H_{12}\) over \(\Q\) by the following equations.
\begin{align*}
 H_4 &= E_4,\qquad K_6=E_6,\\
 H_8 &:= 61(H_4^2-E_8)/(2^{10}\cdot3^2\cdot5^2),\\
 K_{10} &:= 277(H_4K_6-E_{10})/(2^9\cdot3^3\cdot5^2\cdot7),\\
 H_{12} &:= \Bigl(19\cdot2659(441H_4^3+250K_6^2-691E_{12})\\
 &\hspace{8mm}{}-2^9\cdot3^4\cdot5^2\cdot7^2\cdot6791 H_4H_8\Bigr)
 \big/\bigl(2^{11}\cdot3^7\cdot5^3\cdot7^2\cdot73\bigr).
\end{align*}
We use the normalizations in \cite[\S2.5]{Kikuta}.
The algebraic independence of these 5 forms over \(\C\) follows from
the structure theorem of Dern and Krieg~\cite{DK}.
Kikuta and Nagaoka~\cite[Theorems~4.3 and~5.1]{KN} proved that these
forms have integral Fourier coefficients and that, for every prime
\(p\ge5\), they generate the ring of forms with \(p\)-integral Fourier
coefficients.
In particular, these forms generate \(B[1/6]\) over \(\Z[1/6]\). We use their monomials as coordinates in
polynomial calculations below.

\subsection{Siegel modular forms and the restriction map}\label{sec:restriction}

Let \(\mathbb S_2=\{Z\in\mathbb H_2\mid\trn Z=Z\}\). We denote the
Siegel modular group of degree 2 by
\[
 \gamt=\mathrm{Sp}_2(\Z)
      =\{M\in\mathrm{GL}_4(\Z)\mid\trn M J_2M=J_2\}.
\]
The action in Section~\ref{sec:herm-fourier} restricts to an action of
\(\gamt\) on \(\mathbb S_2\). Let \(M_k(\gamt)\) denote the complex vector space of Siegel
modular forms of degree 2 and weight \(k\) on this group. Its elements
are holomorphic functions \(F\colon\mathbb S_2\to\C\) satisfying
\[
 F(M\langle Z\rangle)=\det(CZ+D)^kF(Z)
 \qquad \bigl(M=\left(\smat{A}{B}{C}{D}\right)\in\gamt\bigr).
\]

\begin{samepage}
Every \(F\in M_k(\gamt)\) has a Fourier expansion of the form
\[
 F(Z)=\sum_{0\le T\in\Lambda_2}a_F(T)e^{2\pi i\tr(TZ)},
\]
where \(\Lambda_2\) is the lattice of half-integral symmetric matrices given by
\[
 \Lambda_2=\left\{\mat{m}{r/2}{r/2}{n}\ \middle|\ m,r,n\in\Z\right\}.
\]
We write \(T=(m,r,n)\) for the displayed matrix and \(a_F(m,r,n)\)
for its Fourier coefficient \(a_F(T)\). The condition \(0\le T\) is equivalent
to \(m,n\ge0\) and \(r^2\le4mn\).
\par\end{samepage}
Put \(q_{ij}=e(z_{ij})\) for \(1\le i\le j\le 2\). Since
\(\tr(TZ)=mz_{11}+rz_{12}+nz_{22}\), the expansion becomes
\[
 F(Z)=\sum_{\substack{m,n\ge0,\ r\in\Z\\ r^2\le4mn}}
       a_F(m,r,n)\,q_{11}^m q_{12}^r q_{22}^n .
\]
For a subring \(R\subseteq\C\), the space of forms with Fourier
coefficients in \(R\) is defined by
\[
 M_k(\gamt;R)=
 \bigl\{F\in M_k(\gamt)\mid
   a_F(T)\in R\quad(0\le\forall T\in\Lambda_2)\bigr\}.
\]
We abbreviate \(M_k(\gamt;\Z)\) to \(M_k^{\mathrm S}\).

Since \(\det M=1\) for \(M\in\gamt\), the Hermitian transformation law
on \(\mathbb S_2\) gives the Siegel transformation law above. Thus in
each even weight \(k\) restriction defines the map
\[
 \rho\colon M_k(\gamt(\mathcal O);\Z)\longrightarrow M_k^{\mathrm S},
 \qquad F\longmapsto F|_{\mathbb S_2}.
\]
Setting \(z_{21}=z_{12}\) gives
\[
 \dot q_{11}\mapsto q_{11},\qquad \dot q_{22}\mapsto q_{22},\qquad
 \dot q_{12}\mapsto q_{12},\qquad \ddot q_{12}\mapsto 1.
\]
The Fourier coefficient of \(\rho(F)\) at \(T=(m,r,n)\) is therefore the
sum of the coefficients of \(F\) at all \(H=(m,r,s,n)\) with \(\re H=T\).
Thus
\begin{equation}\label{eq:restriction-coefficients}
 a_{\rho(F)}(m,r,n)
 =\sum_{\substack{s\in\Z\\r^2+s^2\le4mn}}a_F(m,r,s,n).
\end{equation}
In particular, \(\rho\) preserves integral Fourier coefficients.

\subsection{The \(p\)-adic valuation of a Fourier series}\label{sec:valuation}

\begin{samepage}
Fix a prime \(p\). Let \(\ord_p\) be the additive valuation on \(\Q\)
normalized by \(\ord_p(p)=1\), with \(\ord_p(0)=+\infty\).
For a formal Fourier series
\[
 F=\sum_{0\le H\in\Lambda_2(K)}a_F(H)e^{2\pi i\tr(HZ)},
 \qquad a_F(H)\in\Q,
\]
define
\[
 v_p(F):=\inf_{0\le H\in\Lambda_2(K)}\ord_p\bigl(a_F(H)\bigr).
\]
\end{samepage}
Thus \(v_p(0)=+\infty\). For \(F\ne0\), the infimum is a minimum if the
coefficient valuations are bounded below, and is \(-\infty\) otherwise.

Write \(\Z_{(p)}=\{a/b\mid a,b\in\Z,\ p\nmid b\}\) and let \(\F_p\)
denote the field with \(p\) elements. All coefficients of \(F\) belong
to \(\Z_{(p)}\) if and only if \(v_p(F)\ge0\). They belong to \(\Z\)
if and only if this holds for every prime \(p\).
If \(v_p(F),v_p(G)>-\infty\), then
\[
 v_p(F+G)\ge\min\{v_p(F),v_p(G)\},\qquad
 v_p(FG)=v_p(F)+v_p(G).
\]
The inequality follows coefficient by coefficient. The product formula
is \cite[Lemma~2.8]{Kikuta}. Applying the formula for the Fourier
coefficients of \(E_k\) to the definition of \(H_8\) in
Section~\ref{sec:basic-five} gives \(a_{H_8}(1,1,1,1)=1\).
Since all coefficients of \(H_8\) are integers, \(v_p(H_8)=0\) for every prime \(p\).
In particular, if the Fourier coefficients of \(G\) have a common denominator, then
\begin{equation}\label{eq:h8-valuation}
 F=H_8G\quad\Longrightarrow\quad v_p(G)=v_p(F)\qquad\text{for every prime }p .
\end{equation}

\section{Construction, integrality and restrictions of the generators}\label{sec:forms}

\subsection{Defining formulas}\label{sec:generator-definitions}

Starting with the 5 basic forms, we define the remaining generators and
auxiliary forms over \(\Q\). For the generators defined below in weights
divisible by 4 and for \(K_{14},K_{18},K_{22},K_{30},K_{42}\), we use the
constructions in \cite[\S3.1]{Kikuta}.
Our \(K_{10},J_{16},H_{16}\) are denoted there by \(F_{10},H_{16},I_{16}\),
respectively. We prove integrality in Section~\ref{sec:integrality}.

The definitions in weights up to 28 are
\begingroup
\allowdisplaybreaks[2]
\begin{align*}
 I_{12}&:=(H_4^3-K_6^2)/1728+144H_{12},\\
 K_{14}&:=\bigl(H_4K_{10}-K_6H_8\bigr)/6,\\
 J_{16}&:=\bigl(K_6K_{10}-H_4^2H_8\bigr)/6,\\*
 H_{16}&:=\bigl(H_4H_{12}-J_{16}\bigr)/12,\\
 K_{18}&:=\bigl(K_6H_{12}-H_4K_{14}\bigr)/12,\\
 H_{20}&:=\bigl(K_{10}^2-H_4H_8^2-12H_8H_{12}\bigr)/36,\\
 K_{22}&:=\bigl(K_{10}H_{12}-H_8K_{14}\bigr)/6,\\
 H_{24}&:=\bigl(H_{12}^2-H_4H_{20}-4H_8H_{16}\bigr)/24,\\
 K_{26}&:=\bigl(K_{10}H_{16}-H_8K_{18}\bigr)/6,\\
 J_{28}&:=\bigl(K_{10}K_{18}-H_4H_8H_{16}\bigr)/6,\\
 H_{28}&:=\bigl(H_4H_{24}-J_{28}-2H_8^2I_{12}\bigr)/6.
\end{align*}
\endgroup
In weights 30 through 42 we use
\begingroup
\allowdisplaybreaks[2]
\begin{align*}
 K_{30}&:=\bigl(K_6H_{24}-K_{14}H_{16}+2H_8K_{10}I_{12}\bigr)/6,\\
 L_{30}&:=\bigl(K_{10}H_{20}-H_8K_{22}\bigr)/6,\\
 K_{34}&:=\bigl(K_{10}H_{24}-H_8K_{26}\bigr)/6,\\
 H_{36}&:=\bigl(H_{12}H_{24}-H_{16}H_{20}+14H_8H_{28}+6H_8^3H_{12}\bigr)/18,\\
 K_{38}&:=\bigl(K_{10}H_{28}-H_8K_{30}\bigr)/6,\\
 J_{40,1}&:=\bigl(K_{10}K_{30}-H_4H_8H_{28}\bigr)/6,\\
 H_{40}&:=\bigl(H_4H_{36}-J_{40,1}-H_4H_8H_{28}+H_8^3J_{16}+2H_8I_{12}H_{20}\bigr)/4,\\
 K_{42}&:=\bigl(H_{12}K_{30}-K_{14}H_{28}-6H_8I_{12}K_{22}\bigr)/12.
\end{align*}
\endgroup

In weight 46 put
\begin{align*}
 J_{38}&:=H_8^3K_{14}+H_8K_{30}+K_{38},\\*
 K_{46}&:=\bigl(K_{10}H_{36}-H_8J_{38}\bigr)/2.
\end{align*}
In weight 48 the symbol \(H_{10}\) in the earlier formula denotes our \(K_{10}\). Put
\begin{align*}
 J_{40,2}={}&H_{12}H_{28}+2H_{40}+4H_8K_{10}^2H_{12}
                -2H_4H_8^2H_{20}-2H_4H_8^3H_{12}\\
          &+4H_8I_{12}H_{20}+2H_{12}H_8^2I_{12}-H_{16}H_{24}
                -2H_8^3H_{16}+2J_{40,1},\\
 H_{48}:={}&(H_{12}H_{36}-H_{24}^2)/4-H_8J_{40,2}/8.
\end{align*}
Finally, in weight 50 put
\begin{align*}
 J_{42,2}&=K_{10}^3H_{12}+H_8^3K_{18}+K_{42},\\
 J_{42,3}&=-H_{12}K_{30}+K_6H_{36}-H_8H_{16}K_{18},\\
 K_{50}&:=\bigl(K_{10}H_{40}-H_8(3J_{42,2}-2J_{42,3})\bigr)/6.
\end{align*}

\subsection{Integral Fourier coefficients}\label{sec:integrality}

We use the substitution argument of \cite[Lemmas~3.2 and~3.4]{Kikuta}.
There exist integral formal Fourier series
\(S,U,V\) satisfying
\begin{equation}\label{eq:substitution}
 H_4=1+48S,\quad K_6=1+72S+864U,\quad
 K_{10}=6H_{12}+H_4^2H_8-72V.
\end{equation}
For example, substituting for \(K_{10}\) gives
\begin{align*}
 K_{14}&=H_4(H_{12}-12V)+\bigl((H_4^3-K_6)/6\bigr)H_8,\\
 (H_4^3-K_6)/6&=12S+1152S^2+18432S^3-144U.
\end{align*}
Thus in characteristic 0 we have the integral polynomial identity
\[
 K_{14}=(1+48S)(H_{12}-12V)+(12S+1152S^2+18432S^3-144U)H_8 .
\]
For each subsequent definition \(G=N/d\), with \(d\) a positive integer, expand \(N\) using
\eqref{eq:substitution} and the polynomials already obtained in lower
weights or earlier in the same weight. Collecting equal monomials gives
\(N=dP_G\) with \(P_G\in\Z[S,U,V,H_8,H_{12}]\) in every case.
For the denominators in \(K_{46}\) and \(K_{50}\), the same substitutions give
\[
 H_{36}\equiv J_{38}\pmod2,\qquad H_{40}\equiv J_{42,2}\pmod2,\qquad
 H_{40}\equiv J_{42,3}\pmod3.
\]
For example, the common value of \(H_{36}\) and \(J_{38}\) modulo \(2\)
after substitution is
\[
 H_{12}H_8^3+H_8(V+SH_8)^2+H_{12}(V^2+UH_8^2).
\]
As \(K_{10}\equiv H_8\pmod{6}\), the congruence for \(H_{36}\) proves divisibility
by \(2\) in \(K_{46}\). The combination \(3J_{42,2}-2J_{42,3}\) reduces
to \(J_{42,2}\) modulo \(2\) and \(J_{42,3}\) modulo \(3\). This gives
divisibility by \(6\) in \(K_{50}\). The 5 substituted series are
integral, so the calculation establishes integrality of the full
Fourier expansions.

\subsection{Restrictions of the generators}\label{sec:generator-restrictions}

We denote Igusa's generators for the ring of Siegel modular forms of even
weight with integral Fourier coefficients \cite{Igusa} by \(X_k\) and
\(Y_{12}\), using the normalization
of \cite[\S2.4]{Kikuta}. In particular, \(X_4,X_6\) are the
Eisenstein series with constant term 1, and \(X_{10},X_{12}\) are the
cusp forms normalized by \(a_{X_{10}}(1,1,1)=a_{X_{12}}(1,1,1)=1\).
Our normalization satisfies \(Y_{12}=(X_4^3-X_6^2)/1728+144X_{12}\).
The images of the 5 basic forms
are given in \cite[Theorem~2.6]{Kikuta}.
\begin{equation}\label{eq:restriction}
\begin{array}{c|ccccc}
 F&H_4&K_6&H_8&K_{10}&H_{12}\\
 \rho(F)&X_4&X_6&0&6X_{10}&X_{12}
\end{array}
\end{equation}
To see the factor 6 in \(\rho(K_{10})=6X_{10}\), take \(T=(1,1,1)\)
in \eqref{eq:restriction-coefficients}. The formula for the Fourier
coefficients of \(E_k\), applied to the definition of \(K_{10}\), gives
\(a_{K_{10}}(1,1,\pm1,1)=1\) and \(a_{K_{10}}(1,1,0,1)=4\). Hence
\[
 a_{\rho(K_{10})}(1,1,1)
 =\sum_{s=-1}^{1}a_{K_{10}}(1,1,s,1)=1+4+1=6 .
\]

Applying \(\rho\) to the defining identities of Section~\ref{sec:generator-definitions} now gives
\[
\begin{array}{c|cccccc}
 F&I_{12}&K_{14}&J_{16}&H_{16}&K_{18}&H_{20}\\
 \rho(F)&Y_{12}&X_4X_{10}&X_6X_{10}&X_{16}&X_{18}&X_{10}^2
\end{array}
\]
For example, \(6K_{14}=H_4K_{10}-K_6H_8\) gives
\(6\rho(K_{14})=X_4(6X_{10})-X_6\cdot0\). In characteristic 0 this gives
\(\rho(K_{14})=X_4X_{10}\).

The 2 rings in weights divisible by 4 are
\[
 A_{\mathrm H}=\bigoplus_{m\ge0}M_{4m},\qquad
 A_{\mathrm S}=\bigoplus_{m\ge0}M_{4m}^{\mathrm S}.
\]
Every Siegel modular form of weight divisible by 4 with integral Fourier
coefficients lifts to a Hermitian modular form of the same weight with
integral Fourier coefficients. Equivalently,
\begin{equation}\label{eq:surjective}
 \rho(A_{\mathrm H})=A_{\mathrm S}.
\end{equation}
Indeed, the restrictions of the integral forms in
\(\mathcal G_0\setminus\{H_8\}\) give the 10 generators of
\cite[Theorem~5.5]{Kurihara} in our normalization.
Our \(Y_{12}\) differs from Kurihara's by \(144X_{12}\), and \(X_{12}\)
is itself one of the generators.

\section{Generation in every even weight}\label{sec:generation}

\subsection{Weights divisible by 4}\label{sec:mult4}

We first determine which of the forms constructed in
Section~\ref{sec:forms} are needed to generate
\(A_{\mathrm H}=\bigoplus_{m\ge0}M_{4m}\), the ring of symmetric Hermitian
modular forms of weight divisible by 4 with integral Fourier coefficients.
This is the coefficient ring for the module of forms of weight
\(k\equiv2\pmod4\) in Theorem~\ref{thm:main}.

\begin{samepage}
Kurihara~\cite{Kurihara} proved the following result.
We restate it in our notation.

\begin{Prop}\label{prop:mult4}
The 11 forms in \(\mathcal G_0\) generate the ring of Hermitian modular forms
in weights divisible by 4, that is, \(A_{\mathrm H}=\Z[\mathcal G_0]\).
\end{Prop}
\end{samepage}

\begin{proof}[Sketch of proof]
Of Igusa's 14 even-weight generators \cite{Igusa}, 5 have weight \(2\bmod4\).
In a monomial of weight divisible by 4 these occur in pairs, so their
15 pairwise products and the remaining 9 generators give \(9+15=24\)
generators for \(A_{\mathrm S}\). Kurihara used relations among these
forms to reduce the list to the following 10.
\[
 X_4,\ X_{12},\ Y_{12},\ X_{16},\ X_{10}^2,\ X_{24},\
 X_{28},\ X_{36},\ X_{40},\ X_{48}.
\]

These 10 Siegel modular forms are the restrictions, in the displayed order,
of the forms in \(\mathcal G_0\) other than \(H_8\).
Thus, for \(F\in M_k\) with \(4\mid k\), we can lift a polynomial
expression for \(\rho(F)\) to a form \(P\in\Z[\mathcal G_0]\) of weight
\(k\) with \(\rho(P)=\rho(F)\). Over \(\Q\) the kernel of \(\rho\) is
generated by \(H_8\), so \(F-P=H_8G\) with \(G\in M_{k-8}\otimes\Q\).
The coefficients of \(G\) therefore have a common denominator. By \eqref{eq:h8-valuation},
\(v_p(G)=v_p(F-P)\ge0\) for every prime \(p\). Hence \(G\in M_{k-8}\).
For \(k<8\) the difference vanishes, and induction on \(k\) proves
the assertion.
\end{proof}

\subsection{Reducing the restriction to 3 variables}\label{sec:restriction-reduction}

Let \(F\in M_k\) with \(k\equiv2\pmod4\). We first subtract integral
lifts until the restriction has the form \(X_{10}P(Y_{12},X_{36},X_{48})\).
We then prove the divisibility conditions on its coefficients and use
them to construct a lift in \(\sum_{g\in\mathcal G_2}A_{\mathrm H}g\).
The remaining difference has the form \(H_8G\) with \(G\in M_{k-8}\),
so induction lowers the weight by 8.
Our goal is
\[
 M_k=\sum_{g\in\mathcal G_2}M_{k-\operatorname{wt}(g)}\,g,
 \qquad k\equiv2\pmod4,
\]
where \(\operatorname{wt}(g)\) denotes the weight of \(g\) (its subscript), and \(M_j=0\) for \(j<0\).

Proposition~\ref{prop:mult4} and its proof give
\begin{equation}\label{eq:known-restriction}
 \begin{aligned}
 A_{\mathrm H}&=\Z[\mathcal G_0],\qquad \rho(A_{\mathrm H})=A_{\mathrm S},\\
 A_{\mathrm S}&=\Z[X_4,X_{12},Y_{12},X_{16},X_{10}^2,X_{24},X_{28},X_{36},X_{40},X_{48}].
 \end{aligned}
\end{equation}
Thus, for every \(C\in M_{4m}^{\mathrm S}\), we may choose
\(\widetilde C\in M_{4m}\) with \(\rho(\widetilde C)=C\). The defining formulas in Section~\ref{sec:forms}
also give the following restrictions. At this point we use only
these identities; generation by the displayed restrictions will
follow after the coefficient calculation below.
\begin{center}
\begin{tabular}{rcc@{\qquad}rcc}
\hline
Weight & Restriction & Lift & Weight & Restriction & Lift\\
\hline
6&\(X_6\)&\(K_6\)&30&\(X_{10}^3\)&\(L_{30}\)\\
10&\(6X_{10}\)&\(K_{10}\)&34&\(X_{10}X_{24}\)&\(K_{34}\)\\
14&\(X_4X_{10}\)&\(K_{14}\)&38&\(X_{10}X_{28}\)&\(K_{38}\)\\
18&\(X_{18}\)&\(K_{18}\)&42&\(X_{42}\)&\(K_{42}\)\\
22&\(X_{10}X_{12}\)&\(K_{22}\)&46&\(3X_{10}X_{36}\)&\(K_{46}\)\\
26&\(X_{10}X_{16}\)&\(K_{26}\)&50&\(X_{10}X_{40}\)&\(K_{50}\)\\
30&\(X_{30}\)&\(K_{30}\)&&&\\
\hline
\end{tabular}
\end{center}

Let \(F\in M_k\) with \(k\equiv2\pmod4\).
In a monomial of weight \(k\) in Igusa's generators, take out 1 factor
\(X_d\) of weight \(d\equiv2\pmod4\). The remaining product has
weight \(k-d\equiv0\pmod4\). Grouping terms and applying
\eqref{eq:surjective} in weight \(k-d\) gives
\begin{equation}\label{eq:siegel-five}
 \begin{aligned}
 \rho(F)&=C_6X_6+C_{10}X_{10}+C_{18}X_{18}+C_{30}X_{30}+C_{42}X_{42},\\
 C_d&\in M_{k-d}^{\mathrm S}=\rho(M_{k-d})\subset A_{\mathrm S}
       \qquad(d\le k).
 \end{aligned}
\end{equation}
Terms with \(d>k\) are zero.

For \(d\in\{6,18,30,42\}\), choose \(\widetilde C_d\in M_{k-d}\)
with \(\rho(\widetilde C_d)=C_d\), using \eqref{eq:known-restriction}.
Since \(\rho(K_d)=X_d\) for these 4 indices, we have
\[
 \rho\!\left(F-\sum_{d\in\{6,18,30,42\}}\widetilde C_dK_d\right)
 =C_{10}X_{10}.
\]

The table gives an integral lift of \(6X_{10}\). To handle
\(C_{10}X_{10}\), we look for factors of \(C_{10}\) whose product with
\(X_{10}\) already has a lift in the table.
Fix a homogeneous polynomial expression for \(C_{10}\) in the 10
generators in \eqref{eq:known-restriction}. Collect the terms containing
at least one of the 7 factors
\[
 X_4,\quad X_{12},\quad X_{16},\quad X_{10}^2,\quad
 X_{24},\quad X_{28},\quad X_{40}.
\]
Multiply the sum of these terms, with their original coefficients, by
\(X_{10}\) and denote the result by \(R\in M_k^{\mathrm S}\). For each term
choose one listed factor and lift its product
with \(X_{10}\) by the table. The other factors lift by
\eqref{eq:known-restriction}. Let \(\widetilde R\in M_k\) be the sum
of these lifts with their original integer coefficients, so that
\(\rho(\widetilde R)=R\). For example, if \(X_{10}X_{24}Y_{12}^2\)
has coefficient \(c\in\Z\) in \(R\), we include \(cK_{34}I_{12}^2\)
in \(\widetilde R\).

Subtract from \(F\) both the 4 products \(\widetilde C_dK_d\) and
the form \(\widetilde R\). This defines the Hermitian modular form
\[
 F_1:=F-\bigl(\widetilde C_6K_6+\widetilde C_{18}K_{18}
              +\widetilde C_{30}K_{30}+\widetilde C_{42}K_{42}\bigr)
        -\widetilde R\in M_k .
\]
The terms in \(C_{10}\) containing none of the 7 listed factors
involve only \(Y_{12},X_{36},X_{48}\). Write their sum as
\(P(Y_{12},X_{36},X_{48})\). Taking the restriction gives
\begin{equation}\label{eq:remainder}
 \begin{aligned}
 \rho(F_1)&=C_{10}X_{10}-R=X_{10}P(Y_{12},X_{36},X_{48}),\\
 P&=\sum_{i+3j+4\ell=n}c_{ij\ell}Y_{12}^iX_{36}^jX_{48}^{\ell},
       \qquad c_{ij\ell}\in\Z .
 \end{aligned}
\end{equation}
Here \(n=(k-10)/12\), and the sum is over nonnegative integers \(i,j,\ell\).
If this sum is empty, the lifts already subtracted have the same
restriction as \(F\). Otherwise \(n\) is a nonnegative integer, and
it remains to lift \(X_{10}P\).

\subsection{The divisibility conditions on the remaining coefficients}\label{sec:coefficients}

The coefficients in \eqref{eq:remainder} need not be arbitrary:
the remainder is the restriction of the integral form \(F_1\).
We now use this fact to prove the exact conditions needed for the lifts.

\begin{Lem}\label{lem:remainder-coefficients}
If an integral form \(F_1\) satisfies \eqref{eq:remainder}, then
\[
 2\mid c_{i0\ell}\quad\text{for every }i,\ell,
 \qquad
 3\mid c_{ij0}\quad\text{for every }i,j.
\]
In particular, \(6\mid c_{i00}\).
\end{Lem}

\begin{proof}
If all coefficients of \(P\) are zero, the assertion is immediate. Otherwise \(k=10+12n\) with \(n\ge0\).
Over \(\Q\), the 5 basic forms freely generate \(B\otimes\Q\), and
the kernel of \(\rho\) is \((H_8)\), as in the proof of
Proposition~\ref{prop:mult4}. Since \(\rho(K_{10})=6X_{10}\),
\eqref{eq:remainder} gives
\begin{equation}\label{eq:remainder-kernel}
 6F_1-K_{10}P(I_{12},H_{36},H_{48})=H_8G_1,
 \qquad G_1\in\Q[H_4,K_6,H_8,K_{10},H_{12}].
\end{equation}
Multiplying by \(K_6\) gives
\[
 6K_6F_1=K_6K_{10}P(I_{12},H_{36},H_{48})+K_6H_8G_1.
\]
The weight of \(K_6F_1\) is divisible by 4, so Proposition~\ref{prop:mult4}
expresses it as an integral polynomial in \(\mathcal G_0\).
To compare coefficients in this last identity, first write both sides
as rational polynomials in the 5 basic forms and subtract the right-hand side. The resulting polynomial
\(Q\) vanishes on \((H_4,K_6,H_8,K_{10},H_{12})\). Since these forms are
algebraically independent, \(Q\) is the zero polynomial. It remains zero
under every polynomial substitution.

Let \(T_6,T_{10},T_{12}\) be indeterminates. In the polynomial identity
we make the substitution
\begin{equation}\label{eq:coefficient-substitution}
 (H_4,K_6,H_8,K_{10},H_{12})
 \longmapsto (12,12T_6,T_6T_{10},12T_{10},12T_{12}).
\end{equation}
In particular \(Q(12,12T_6,T_6T_{10},12T_{10},12T_{12})=0\).
The choice \(12^3=1728\) simplifies \(I_{12}\), and the images of
\(K_6K_{10}\) and \(H_4^2H_8\) agree, so \(J_{16}\) maps to 0.

Write \(\operatorname{coeff}_{i,j}(F)\in\Q[T_{12}]\) for the coefficient
of \(T_6^iT_{10}^j\) after substitution. In
\(K_6\) times \eqref{eq:remainder-kernel}, the term \(K_6H_8G_1\)
has a factor \(T_6^2\). Its \(T_6T_{10}\)-coefficient is zero, whereas
\(K_6K_{10}\) maps to \(144T_6T_{10}\). Thus
\begin{equation}\label{eq:coefficient-identity}
 6\operatorname{coeff}_{1,1}(K_6F_1)
 =144\operatorname{coeff}_{0,0}\!\bigl(P(I_{12},H_{36},H_{48})\bigr).
\end{equation}

Substitution in the defining identities gives the following entries.
The first row covers \(j=0,1,2,3\).
\begin{center}
\begin{tabular}{ccc}
\hline
\(F\) & \(\operatorname{coeff}_{0,0}(F)\) & \(\operatorname{coeff}_{1,1}(F)\)\\
\hline
\(H_{4+12j}\)&\(12T_{12}^{j}\)&0\\
\(H_8\)&0&1\\
\(H_{20}\)&0&\(-4T_{12}\)\\
\(H_{12}\)&\(12T_{12}\)&0\\
\(I_{12}\)&\(1+1728T_{12}\)&0\\
\(H_{24}\)&\(6T_{12}^2\)&0\\
\(H_{36}\)&\(4T_{12}^3\)&\(12T_{12}^2\)\\
\(H_{48}\)&\(3T_{12}^4\)&\(24T_{12}^3\)\\
\hline
\end{tabular}
\end{center}
For example, \(12H_{16}=H_4H_{12}-J_{16}\) gives the image
\(12T_{12}\) of \(H_{16}\). The identity for \(36H_{20}\) gives
\[
 36H_{20}\longmapsto
 144T_{10}^2-12T_6^2T_{10}^2-144T_6T_{10}T_{12},
\]
whose \(T_6T_{10}\)-coefficient is \(36(-4T_{12})\).

For each generator of \(A_{\mathrm H}\), the exponents of \(K_6,K_{10}\)
in each basic monomial have even sum, by weight modulo 4. After
\eqref{eq:coefficient-substitution} the same holds for \(T_6,T_{10}\).
Thus the coefficient of \(T_6T_{10}\) in a product comes from
the last column in exactly 1 factor and the constant column in all
other factors.

Consider a product of original weight \(4\bmod12\). All nonzero
constants outside the first row have weight \(0\bmod12\).
\begin{enumerate}[1.]
\item If the last-column factor is \(H_{36}\) or \(H_{48}\), its entry
is a multiple of 12. At least 1 constant comes from the first row,
giving another factor of 12.
\item If it is \(H_8\) or \(H_{20}\), of weight \(8\bmod12\), the
constants must have total weight \(8\bmod12\). At least 2 come from
the first row, giving a factor of \(12^2\).
\end{enumerate}
Since \(K_6F_1\in A_{\mathrm H}\) has weight \(16+12n\), we have
\(\operatorname{coeff}_{1,1}(K_6F_1)=144R\) for some
\(R\in\Z[T_{12}]\). Substituting in \eqref{eq:coefficient-identity}
and comparing the two multiples of 144 in characteristic 0 gives
\begin{equation}\label{eq:remainder-congruence}
 P(1+1728T_{12},4T_{12}^3,3T_{12}^4)=6R .
\end{equation}
We can now reduce the coefficients in \(\Z[T_{12}]\). Modulo 2 and
modulo 3, respectively, we obtain
\[
 \sum_{i+4\ell=n}c_{i0\ell}T_{12}^{4\ell}\equiv0\pmod2,\qquad
 \sum_{i+3j=n}c_{ij0}T_{12}^{3j}\equiv0\pmod3 .
\]
In the first sum each exponent determines \(\ell\) and then \(i\).
In the second it determines \(j\) and then \(i\). Comparison of
coefficients proves the 2 divisibility assertions.
\end{proof}

\subsection{Constructing the remaining lifts}\label{sec:remaining-lifts}

The conditions just proved give 4 cases. The following table lists
the basic product to lift in each case.
\begin{center}
\begin{tabular}{cccc}
\hline
Exponents & Divisibility & Basic product & Lift\\
\hline
\(j=\ell=0\)&\(6\mid c_{i00}\)&\(6X_{10}\)&\(K_{10}\)\\
\(j>0,\ \ell=0\)&\(3\mid c_{ij0}\)&\(3X_{10}X_{36}\)&\(K_{46}\)\\
\(j=0,\ \ell>0\)&\(2\mid c_{i0\ell}\)&\(2X_{10}X_{48}\)&\(J_{58}\)\\
\(j>0,\ \ell>0\)&\(c_{ij\ell}\in\Z\)&\(X_{10}X_{36}X_{48}\)&\(J_{94}\)\\
\hline
\end{tabular}
\end{center}
The 2 auxiliary forms here are the following explicit combinations,
so they introduce no additional generators:
\begin{align*}
 J_{58}&=K_{10}H_{48}-K_{22}H_{36}+K_{34}H_{24},\\
 J_{94}&=K_{22}H_{36}^{2}-K_{34}H_{24}H_{36}-K_{46}H_{48}.
\end{align*}
Indeed, Igusa's generator \(X_{48}\) is defined by
\(4X_{48}=X_{12}X_{36}-X_{24}^2\). Hence
\begin{align*}
 \rho(J_{58})&=(6-4)X_{10}X_{48}=2X_{10}X_{48},\\
 \rho(J_{94})&=(4-3)X_{10}X_{36}X_{48}=X_{10}X_{36}X_{48}.
\end{align*}
For a term \(c_{ij\ell}X_{10}Y_{12}^iX_{36}^jX_{48}^{\ell}\) in
\eqref{eq:remainder}, take the following lift. Its scalar coefficient
is integral by Lemma~\ref{lem:remainder-coefficients}.
\[
 \begin{cases}
 (c_{i00}/6) I_{12}^iK_{10},&j=\ell=0,\\[2pt]
 (c_{ij0}/3) I_{12}^iH_{36}^{j-1}K_{46},&j>0,\ \ell=0,\\[2pt]
 (c_{i0\ell}/2) I_{12}^iH_{48}^{\ell-1}J_{58},&j=0,\ \ell>0,\\[2pt]
 c_{ij\ell} I_{12}^iH_{36}^{j-1}H_{48}^{\ell-1}J_{94},&j>0,\ \ell>0.
 \end{cases}
\]
For example, in weight 46 the remainder has the form
\[
 c_{300}X_{10}Y_{12}^3+c_{010}X_{10}X_{36},\qquad
 6\mid c_{300},\quad 3\mid c_{010},
\]
and its lift is \((c_{300}/6)I_{12}^3K_{10}+(c_{010}/3)K_{46}\).

Adding these lifts to those subtracted in
Section~\ref{sec:restriction-reduction} gives a homogeneous form
\begin{equation}\label{eq:chosen-lift}
 F_0\in\sum_{g\in\mathcal G_2}A_{\mathrm H}g,
 \qquad \rho(F_0)=\rho(F).
\end{equation}
This completes the construction of a lift with the same restriction.

\subsection{Completing generation by induction}\label{sec:final-induction}

Let \(F\in M_k\), with \(k\equiv2\pmod4\), and take \(F_0\)
as in \eqref{eq:chosen-lift}. Since the rational kernel of \(\rho\)
is \((H_8)\),
\[
 \rho(F-F_0)=0\quad\Longrightarrow\quad
 F-F_0=H_8G,\qquad G\in M_{k-8}\otimes\Q.
\]
As in Proposition~\ref{prop:mult4}, \(G\) has bounded Fourier
denominators. By \eqref{eq:h8-valuation},
\[
 v_p(G)=v_p(F-F_0)\ge0\qquad\text{for every prime }p.
\]
Thus \(G\in M_{k-8}\). For \(k<8\) the rational kernel in
weight \(k\) is zero. Induction on \(k\), using
\(F=F_0+H_8G\) and \(H_8\in A_{\mathrm H}\), proves
\[
 \bigoplus_{m\ge0}M_{4m+2}
 =\sum_{g\in\mathcal G_2}A_{\mathrm H}g.
\]
Together with Proposition~\ref{prop:mult4}, this proves the generation
assertion in Theorem~\ref{thm:main}.

In weight 10 this result explains why lifting \(C_{10}X_{10}\)
required a separate argument in Section~\ref{sec:restriction-reduction}. The only generator monomials of this weight are \(H_4K_6\) and \(K_{10}\), so \eqref{eq:restriction} gives
\[
 \rho(M_{10})=\Z X_4X_6+6\Z X_{10}\ \subsetneq\ M_{10}^{\mathrm S}.
\]
Indeed, an integral lift of \(X_{10}\) would give
\(X_{10}=aX_4X_6+6bX_{10}\) with \(a,b\in\Z\). Comparing
constant terms gives \(a=0\), hence \(6b=1\).

\section{Minimality}\label{sec:minimality}

Section~\ref{sec:generation} gives 24 generators. To prove that every
generating set has at least 24 elements, we count what remains in each
weight after removing products of lower weights and reducing modulo 2.
Let \(D_k\subseteq M_k\) be the subgroup spanned by products of forms
of smaller positive weights whose weights sum to \(k\). For \(k>0\)
and a prime \(p\), let \(Q_{k,p}\) be the \(\F_p\)-space of classes in
\(M_k\) under the equivalence
\[
 F\sim G\quad\Longleftrightarrow\quad F-G\in pM_k+D_k.
\]
In \(\bigoplus_{k>0}Q_{k,p}\) every product of at least 2 positive-weight
elements has zero image. Therefore the images of any algebra generators
of \(B\), after subtracting their components of weight 0, span this
direct sum. The same argument applies to inhomogeneous generators.

\subsection{An explicit matrix calculation in weight 14}\label{sec:weight14}

We give the complete matrix calculation in 1 weight before applying
the same procedure to the remaining weights.

Let \(F\in M_{14}\). The 5 basic forms have exactly 3 monomials of weight 14,
namely \(H_4^2K_6\), \(H_8K_6\) and \(H_4K_{10}\). These give a basis
of \(M_{14}\otimes\Q\). The identity \(6K_{14}=H_4K_{10}-H_8K_6\) allows us to replace
\(H_4K_{10}\) by \(K_{14}\) in this basis. Thus there are unique
\(c_1,c_2,c_3\in\Q\) such that
\[
 F=c_1H_4^2K_6+c_2H_8K_6+c_3K_{14}.
\]
We want to show that \(c_1,c_2,c_3\) are integers. Multiplication by
\(K_6\) sends \(F\) to weight 20, where Proposition~\ref{prop:mult4}
already gives us an integral generating set. We shall identify \(c_1,c_2,c_3\) with coefficients in
the expansion of \(K_6F\) in a \(\Z\)-basis of \(M_{20}\).

To obtain a convenient \(\Z\)-basis of \(M_{20}\), put
\(J_{12}=I_{12}-144H_{12}=(H_4^3-K_6^2)/1728\).
We may replace \(I_{12}\) by \(J_{12}\) among the generators of \(A_{\mathrm H}\),
since \(H_{12}\) is also a generator. Proposition~\ref{prop:mult4}
shows that the following 9 monomials span \(M_{20}\) over \(\Z\).
There are also 9 monomials of weight 20 in the 5 algebraically
independent basic forms, so \(\dim_{\Q}(M_{20}\otimes\Q)=9\).
The displayed monomials therefore form a \(\Z\)-basis of \(M_{20}\).
\[
\begin{gathered}
 H_4^5,\quad H_4^3H_8,\quad H_4^2H_{12},\quad H_4^2J_{12},
 \quad H_4H_{16},\quad H_8J_{12},\\
 H_4H_8^2,\quad H_8H_{12},\quad H_{20}.
\end{gathered}
\]

In \eqref{eq:weight20}, each row expands the product on the left in
the 9 basis monomials listed across the top. For example, row 1 is the identity
\[
 H_4^2K_6^2=H_4^2(H_4^3-1728J_{12})=H_4^5-1728H_4^2J_{12}.
\]
The entries are \(1\) under \(H_4^5\), \(-1728\) under \(H_4^2J_{12}\),
and \(0\) elsewhere. The remaining rows follow from
\[
 K_6^2=H_4^3-1728J_{12},\qquad
 K_6K_{10}=H_4^2H_8+6H_4H_{12}-72H_{16}.
\]
For the last row, substitution gives
\begin{align*}
 6K_6K_{14}
 &=H_4(K_6K_{10})-H_8K_6^2\\
 &=6H_4^2H_{12}-72H_4H_{16}+1728H_8J_{12}.
\end{align*}
Both sides are 6 times an integral form. Comparing these forms in
characteristic 0 gives row 3.
\begingroup
\setlength{\arraycolsep}{2.5pt}
\begin{equation}\label{eq:weight20}
 \begin{array}{@{}l@{\quad}rrrrrrrrr@{}}
 \text{Product}&H_4^5&H_4^3H_8&H_4^2H_{12}&H_4^2J_{12}&H_4H_{16}
     &H_8J_{12}&H_4H_8^2&H_8H_{12}&H_{20}\\ \hline
 H_4^2K_6^2    &1&0&0&-1728&0&0&0&0&0\\
 H_8K_6^2      &0&1&0&0&0&-1728&0&0&0\\
 K_6K_{14}     &0&0&1&0&-12&288&0&0&0
 \end{array}
\end{equation}
\endgroup
Write \(R_i\) for row \(i\) of \eqref{eq:weight20}. The coefficient row
of \(K_6F\) in the integral basis of \(M_{20}\) is
\[
 c_1R_1+c_2R_2+c_3R_3
   =(c_1,c_2,c_3,-1728c_1,-12c_3,-1728c_2+288c_3,0,0,0)
   \in\Z^9.
\]
The columns headed \(H_4^5,H_4^3H_8,H_4^2H_{12}\) give
\(c_1,c_2,c_3\), respectively. Since \(K_6F\in M_{20}\) and
the columns index a \(\Z\)-basis, these 3 coefficients are integers. Hence
\[
 M_{14}=\Z H_4^2K_6\oplus\Z H_8K_6\oplus\Z K_{14}.
\]

The monomials
of weight 14 in the lower-weight generators are
\(H_4^2K_6,\allowbreak H_8K_6,\allowbreak H_4K_{10}\). In the basis just obtained, the identity
\[H_4K_{10}=H_8K_6+6K_{14}\]
shows that their respective coefficients at \(K_{14}\) are 0, 0 and 6.
Thus the coefficient of \(K_{14}\) in every element of \(D_{14}\) is a
multiple of 6, and \(6K_{14}\in D_{14}\). The quotient of \(M_{14}\) by these
products is therefore cyclic of order 6, generated by the class
of \(K_{14}\).
In particular, \(\dim_{\F_2}Q_{14,2}=1\).

\subsection{The 2 generators of weight 30}

By Section~\ref{sec:generation},
\(M_{30}=D_{30}+\Z K_{30}+\Z L_{30}\). We determine for which integers
\(m,n\) the form \(mK_{30}+nL_{30}\) belongs to \(D_{30}\).
Start with the defining identities
\[
\begin{aligned}
 6K_{30}&=K_6H_{24}-K_{14}H_{16}+2H_8K_{10}I_{12},\\
 6L_{30}&=K_{10}H_{20}-H_8K_{22}.
\end{aligned}
\]
Since \(6K_{14}=H_4K_{10}-K_6H_8\), all terms on the right apart
from \(K_6H_{24}\) and \(K_{10}H_{20}\) contain \(H_4\) or \(H_8\)
when expanded in the 5 basic forms. We therefore keep the monomials
containing neither factor. The 5 forms are algebraically independent,
so this amounts to setting \(H_4=H_8=0\) in the unique polynomial
expansion. Under this substitution the definitions give
\(H_{24}\mapsto H_{12}^2/24\) and \(H_{20}\mapsto K_{10}^2/36\). Hence
\[
\begin{aligned}
 K_{30}&\longmapsto (K_6/6)(H_{12}^2/24)=K_6H_{12}^2/144,\\
 L_{30}&\longmapsto (K_{10}/6)(K_{10}^2/36)=K_{10}^3/216.
\end{aligned}
\]
Thus in \(mK_{30}+nL_{30}\) the coefficient of \(K_6H_{12}^2\)
is \(m/144\), and that of \(K_{10}^3\) is \(n/216\).

We compare these coefficients with those of the products spanning
\(D_{30}\). As in Section~\ref{sec:weight14}, we use
\(J_{12}=I_{12}-144H_{12}\) in place of \(I_{12}\). Of the 49 monomials
of weight 30 in the lower-weight generators, the following 5 have a
nonzero coefficient at \(K_6H_{12}^2\) or \(K_{10}^3\).
\[
 \begin{array}{@{}l@{\quad}rr@{}}
 \text{Product}&K_6H_{12}^2&K_{10}^3\\ \hline
 K_6H_{12}^2&1&0\\
 H_{12}K_{18}&1/12&0\\
 K_6H_{24}&1/24&0\\
 K_{10}^3&0&1\\
 K_{10}H_{20}&0&1/36
 \end{array}
\]
Thus for every element of \(D_{30}\), the coefficient of \(K_6H_{12}^2\)
lies in \((1/24)\Z\), and that of \(K_{10}^3\) lies in \((1/36)\Z\).
Consequently \(mK_{30}+nL_{30}\in D_{30}\) forces
\[
 m/144\in(1/24)\Z,\qquad n/216\in(1/36)\Z,
\]
or equivalently \(6\mid m\) and \(6\mid n\).
These conditions are also sufficient, since the defining identities
above express \(6K_{30}\) and \(6L_{30}\) as sums of lower-weight
products with integer coefficients. It follows that
\(M_{30}/D_{30}\simeq(\Z/6\Z)^2\), with generators given by the
classes of \(K_{30},L_{30}\). In particular, \(\dim_{\F_2}Q_{30,2}=2\).

\subsection{The lower bound in all weights}\label{sec:weight-ranks}

For each weight \(k\) occurring in the generating set, expand all
weight-\(k\) monomials in \(\mathcal G_0\cup\mathcal G_2\) in the
5 basic forms. Their integer span is \(M_k\) by
Section~\ref{sec:generation}. Express the products of length at least
2 in an integer basis of this span. Each product gives a row of integer
coefficients. Reduce these entries modulo \(p\). Write \(d_k=\rank M_k\) and let \(r_{k,p}\) be the rank over
\(\F_p\) of these rows. Then
\[
 \dim_{\F_p}Q_{k,p}=d_k-r_{k,p}.
\]
For \(k=14\), the first 2 basis elements are lower-weight products, and the remaining product satisfies \(H_4K_{10}\equiv H_8K_6\pmod p\) for \(p=2,3\). Thus \(r_{14,2}=r_{14,3}=2\). The same calculation in the other generator
weights gives the following ranks. We include \(p=3\) for
Corollary~\ref{cor:graded-small-primes}.

\begin{center}
\begin{tabular}{@{}rrrr@{\qquad\qquad}rrrr@{}}
\hline
\(k\)&\(d_k\)&\(r_{k,2}\)&\(r_{k,3}\)&\(k\)&\(d_k\)&\(r_{k,2}\)&\(r_{k,3}\)\\
\hline
4 & 1 & 0 & 0 & 26 & 13 & 12 & 12 \\
6 & 1 & 0 & 0 & 28 & 19 & 18 & 18 \\
8 & 2 & 1 & 1 & 30 & 20 & 18 & 18 \\
10 & 2 & 1 & 1 & 34 & 27 & 26 & 26 \\
12 & 4 & 2 & 2 & 36 & 36 & 35 & 35 \\
14 & 3 & 2 & 2 & 38 & 36 & 35 & 35 \\
16 & 6 & 5 & 5 & 40 & 47 & 46 & 47 \\
18 & 6 & 5 & 5 & 42 & 49 & 48 & 49 \\
20 & 9 & 8 & 8 & 46 & 63 & 62 & 63 \\
22 & 9 & 8 & 8 & 48 & 78 & 77 & 78 \\
24 & 14 & 13 & 13 & 50 & 80 & 79 & 80 \\
\hline
\end{tabular}
\end{center}

All monomials of other positive weights have length at least 2.
For \(p=2\), weights 12 and 30 each contribute 2 dimensions, and the
other 20 weights in the table contribute 1 each. Hence
\[
 \sum_{k>0}\dim_{\F_2}Q_{k,2}=24 .
\]
Every generating set therefore has at least 24 elements.
The 24 forms in \(\mathcal G_0\cup\mathcal G_2\) attain this bound.
This completes the proof of Theorem~\ref{thm:main}.
\hfill\qedsymbol

\begin{Rem}[Integral lifts in weight 46]\label{rem:weight46}
In \cite[Remark~3.8]{Kikuta} we predicted that \(X_{10}X_{36}\) has
no Hermitian lift with integral Fourier coefficients and observed that
\(3X_{10}X_{36}\) does have such a lift. We can now determine exactly
which multiples lift. Namely,
\[
 mX_{10}X_{36}\in\rho(M_{46})
 \quad\Longleftrightarrow\quad 3\mid m
 \qquad(m\in\Z).
\]
Indeed, if \(mX_{10}X_{36}\) lifts, Lemma~\ref{lem:remainder-coefficients}
with \(P=mX_{36}\) gives \(3\mid m\). Conversely, if \(3\mid m\), then \((m/3)K_{46}\) is an integral lift since
\(\rho(K_{46})=3X_{10}X_{36}\).
Taking \(m=1\) proves the prediction.

For comparison, the least positive multipliers for \(X_{10}X_{40}\)
and \(X_{10}X_{48}\) are 1 and 2, respectively. Indeed, \(K_{50}\)
lifts the first product. For the second, the same lemma with
\(P=mX_{48}\) gives \(2\mid m\), and \(J_{58}\) in
Section~\ref{sec:remaining-lifts} lifts \(2X_{10}X_{48}\).
\end{Rem}

\section{Sturm bounds and reduction modulo 2 and 3}\label{sec:sturm-and-reduction}

\subsection{The diagonal Sturm bound in even weights}\label{sec:sturm}

For completeness, we record the diagonal Sturm bound in even weights.
The proof is independent of Theorem~\ref{thm:main}.
The diagonal bound \(\lfloor k/8\rfloor\) is known
for \(p\ge5\) in all even weights and for every prime when \(4\mid k\)
\cite[Theorems~2.10 and~3.9]{Kikuta}. In the latter case the proof for
\(p=2,3\) uses integral lifts. Frobenius gives the remaining
case \(p=2,3\), \(k\equiv2\pmod4\). The products below show sharpness at
every prime. For the corresponding bound and examples in characteristic~0,
see \cite[Theorem~1.3 and Example~5.7]{NST}.

\begin{Prop}\label{prop:sturm}
Let \(p\) be a prime, let \(k\ge0\) be even, and put \(d=\lfloor k/8\rfloor\). Let \(F\) be a symmetric Hermitian modular
form on \(\gamt(\mathcal O)\) of weight \(k\) and character
\(\det^{k/2}\) with Fourier coefficients in \(\Z_{(p)}\). If
\[
 a_F(m,r,s,n)\equiv0\pmod p
 \quad\text{for all }\quad
 0\le m,n\le k/8,\qquad r^2+s^2\le4mn,
\]
then \(F\equiv0\pmod p\).

For each even weight \(k\ne2\), the bound \(d\) is sharp for every
prime \(p\). More precisely, there is an integral form \(F_k\in M_k\),
independent of \(p\), such that
\[
 a_{F_k}(m,r,s,n)=0\quad\text{if }m<d\text{ or }n<d,
 \qquad a_{F_k}(d,d,d,d)=1 .
\]
In weight 2, \(M_2=0\).
\end{Prop}

\begin{proof}
\textit{The upper bound.}
Weight 0 consists of constants. By the results cited above, it suffices
to consider \(k\equiv2\pmod4\) and \(p=2,3\). Suppose that \(F\)
satisfies the coefficient hypothesis. Every term of \(F\) surviving
modulo \(p\) has \(m\ge d+1\) or \(n\ge d+1\). For a power \(q\) of
\(p\), the Frobenius identity for Fourier expansions gives
\[
 F(Z)^q\equiv \sum_{H\ge0}a_F(H)^q e^{2\pi i\tr(qHZ)}\pmod p.
\]
Thus every surviving term of \(F^q\) has \(m\ge q(d+1)\) or
\(n\ge q(d+1)\). Multiplication by an integral form preserves this
condition since its Fourier indices have nonnegative diagonal entries.

For \(p=2\), take \(F^2\), of weight \(2k\) divisible by 4.
Its coefficients vanish modulo 2 for \(0\le m,n\le 2d+1\).
Since \(k\le8d+6\),
\[
 \lfloor k/4\rfloor\le\lfloor(16d+12)/8\rfloor=2d+1.
\]
Thus the bound in weight \(2k\) is at most \(2d+1\). It gives
\(F^2\equiv0\pmod2\), hence
\(F\equiv0\pmod2\) by the same Frobenius identity.

For \(p=3\), take \(K_6F^9\). Its weight \(9k+6\) is divisible by 4.
We use the ninth power so that even after adding weight 6 the bound
is strictly below the least possible value of \(\max\{m,n\}\) for a
surviving term. Indeed,
\[
 \lfloor(9k+6)/8\rfloor\le\lfloor(72d+60)/8\rfloor=9d+7<9(d+1).
\]
Every coefficient in the required square therefore vanishes modulo 3.
The bound in weight \(9k+6\) gives \(K_6F^9\equiv0\pmod3\).
By \eqref{eq:substitution}, \(K_6\equiv1\pmod3\), so \(F^9\equiv0\pmod3\).
The Frobenius identity then gives \(F\equiv0\pmod3\).

\medskip\noindent\textit{Sharpness.}
Using the formula for the Fourier coefficients of the Hermitian Eisenstein
series \(E_k\) given in \cite[Theorem~2.1]{Kikuta}, we obtain the following
expansions in the variables of Section~\ref{sec:herm-fourier}.
\begin{align*}
 H_8&=\dot q_{11}\dot q_{22}\bigl((\dot q_{12}+\dot q_{12}^{-1}-2)(\ddot q_{12}+\ddot q_{12}^{-1}-2)+\cdots\bigr),\\
 K_{10}&=\dot q_{11}\dot q_{22}\bigl((\dot q_{12}+\dot q_{12}^{-1}+4)(\ddot q_{12}+\ddot q_{12}^{-1}+4)-36+\cdots\bigr).
\end{align*}
The omitted terms have \(m>1\) or \(n>1\) after multiplication by
\(\dot q_{11}\dot q_{22}\). On the axes \(mn=0\) we have \(r=s=0\),
so \eqref{eq:restriction-coefficients} and the restrictions \(0\)
and \(6X_{10}\) give zero coefficients for both forms.
For \(k\ne2\), take the integral form in the following table.
\[
 \begin{array}{c|cccc}
 k & 8d & 8d+2\ (d\ge1) & 8d+4 & 8d+6\\ \hline
 F_k & H_8^d & H_8^{d-1}K_{10} & H_8^dH_4 & H_8^dK_6
 \end{array}
\]
Each product has \(d\) factors from \(H_8,K_{10}\), so its terms have
\(m,n\ge d\). To obtain \((d,d,d,d)\), each of these factors must
contribute a term with \(m=n=1\). In the displayed terms with \(m=n=1\), the largest
exponent in each off-diagonal variable is 1. Thus each factor must
contribute \((1,1,1,1)\), whose coefficient is 1. The remaining
factor, if present, contributes its constant term 1.
Thus \(a_{F_k}(d,d,d,d)=1\) survives modulo every prime.
For \(d\ge1\), all coefficients in the smaller square
\(0\le m,n\le d-1\) vanish. For \(d=0\), the constant term is 1.
By \cite[Theorem~2.6]{Kikuta}, the graded algebra \(B\otimes\Q\) is
generated in weights 4, 6, 8, 10 and 12. There is no monomial of
weight 2 in these generators, so \(M_2=0\).
\end{proof}

\subsection{The ring of Fourier expansions modulo 2 and 3}\label{sec:residue-rings}

For an integral form \(F\), write \(\overline F\) for its Fourier series
with coefficients reduced modulo \(p\). Let \(\widetilde B_p\) be the
image of \(B\) under this reduction, so forms of different weights with
congruent Fourier expansions have the same image. For \(p\ge5\) the
structure of this ring is given in \cite[Theorem~5.2]{KN}.
For related work on Hermitian modular forms over \(\Q(i)\) modulo \(p\), see
Meher and Singh~\cite{MS}.

\begin{Cor}\label{cor:fourier-small-primes}
For \(p=2,3\), the ring of Fourier expansions modulo \(p\) is the
polynomial ring
\[
 \widetilde B_p=
 \F_p[\overline H_8,\overline H_{12},\overline I_{12},\overline H_{16}].
\]
\end{Cor}

\begin{proof}
By \eqref{eq:substitution}, \(H_4\equiv1\pmod{6}\) and \(K_6\equiv1\pmod{6}\). If
\(F\in M_k\) and \(k\equiv2\pmod4\), then \(K_6F\in A_{\mathrm H}\)
and \(K_6F\equiv F\pmod6\). Thus the 11 generators of
Proposition~\ref{prop:mult4} suffice. Of these, \(H_4\) reduces to 1,
and \(H_8,H_{12},I_{12},H_{16}\) are the proposed generators. It remains
to express \(H_{20},H_{24},H_{28},H_{36},H_{40},H_{48}\) in those 4.

We compute all 6 expansions modulo 6. Substituting
\eqref{eq:substitution} in the identities for \(I_{12}\) and \(H_{16}\)
over \(\Z\) gives
\begin{align*}
 I_{12}&=144H_{12}+S^2+64S^3-U-72SU-432U^2,\\
 H_{16}&=K_6V-(2S+72U)H_{12}-(S+12U)H_4^2H_8.
\end{align*}
Hence modulo 6 we have
\[
 U\equiv S^2+4S^3-I_{12}\pmod{6},\qquad
 V\equiv H_{16}+S(H_8+2H_{12})\pmod{6}.
\]
Use these in the integer polynomials \(P_G(S,U,V,H_8,H_{12})\) of
Section~\ref{sec:integrality} and collect terms modulo 6.
For example, the first congruence below also follows at once from
\(H_4H_{20}=H_{12}^2-4H_8H_{16}-24H_{24}\).
In the first 3 weights we obtain
\begin{align*}
 H_{20}&\equiv H_{12}^2+2H_8H_{16}\pmod{6},\\
 H_{24}&\equiv H_{12}H_{16}-2H_8^2 I_{12}\pmod{6},\\
 H_{28}&\equiv H_{16}^2-2H_{12}H_8 I_{12}\pmod{6}.
\end{align*}
In the last 3 weights the congruences modulo 6, grouped by powers
of \(I_{12}\), are
\begingroup
\allowdisplaybreaks[2]
\begin{align*}
 H_{36}\equiv{}&
 (H_{12}+H_8)H_{16}^2-H_{12}H_8^3-2H_{16}^3\\
 &+I_{12}(2H_{12}^3+2H_{12}^2H_8+3H_{12}H_8^2+2H_{16}H_8^2)
       -2I_{12}^2H_8^3\pmod{6},\\[2pt]
 H_{40}\equiv{}&
 H_{16}^3+3H_{16}H_8^3-H_{12}H_8^3\\
 &+I_{12}(2H_{12}^3+H_{12}^2H_8-H_{16}H_8^2)-2I_{12}^2H_8^3\pmod{6},\\[2pt]
 H_{48}\equiv{}&
 H_{12}H_{16}^3+3H_{16}^4+3H_{12}^2H_8^3+2H_{12}H_8^4+2H_{16}^3H_8\\
 &+I_{12}(-H_{12}^4+3H_{12}^3H_8+H_{12}^2H_8^2+
               3H_{12}H_{16}H_8^2+H_{16}H_8^3)\\
 &+I_{12}^2(H_8^4-2H_{12}H_8^3)\pmod{6}.
\end{align*}
\endgroup
Every right-hand side is a polynomial in the stated 4 forms. Reducing
modulo \(p\) now proves generation for both \(p=2\) and \(p=3\).

Fix \(p=2\) or \(3\). To prove algebraic independence over \(\F_p\),
restrict to the Siegel upper half-space.
By \eqref{eq:restriction} and the definitions in Section~\ref{sec:forms},
the restrictions of \(H_8,H_{12},I_{12},H_{16}\) are
\(0,X_{12},Y_{12},X_{16}\), respectively. For either prime,
\(\overline X_{12}=\overline X_{10}\), and
\(\overline X_{10},\overline Y_{12},\overline X_{16}\) are algebraically
independent by \cite[Lemmas~1--2 and Theorem~2]{Nagaoka}. If a nonzero
relation exists, choose a nonzero polynomial \(P(T_1,T_2,T_3,T_4)\)
of least total degree giving such a relation. Restricting this relation and using the algebraic independence of
the 3 reduced Fourier series gives \(P(0,T_2,T_3,T_4)=0\) as a polynomial.
Hence \(P=T_1Q\), and
\[
 \overline H_8\,Q(\overline H_8,\overline H_{12},
                         \overline I_{12},\overline H_{16})=0 .
\]
Since \(a_{H_8}(1,1,1,1)=1\), the first factor is nonzero. The Fourier
series ring over \(\F_p\) is a domain, so \(Q\) is a relation of
smaller degree, a contradiction.
\end{proof}

\subsection{Generators of the graded algebra modulo 2 and 3}\label{sec:graded-reductions}

In the graded algebra \(B_p=B\otimes_{\Z}\F_p\), each form retains its
weight. Sending each class to its reduced Fourier expansion gives a
surjection \(B_p\to\widetilde B_p\). Two forms \(F,G\in M_k\) have the same class precisely when
\(F-G\in pM_k\). For instance, the classes of \(H_4\) and \(K_6\) lie in weights
4 and 6, although both have Fourier expansion 1 modulo \(p\).

\begin{Cor}\label{cor:graded-small-primes}
The minimum numbers of algebra generators of \(B_2\) and \(B_3\) over
their respective prime fields are 24 and 19. Modulo 3, the generators
of Theorem~\ref{thm:main} of weights at most 38 form a minimal set.
\end{Cor}

\begin{proof}
As in Section~\ref{sec:minimality}, the images of any generating set span
\(\bigoplus_{k>0}Q_{k,p}\) after the weight-0 components are removed.
Its weight-\(k\) component has dimension \(d_k-r_{k,p}\).
The ranks in Section~\ref{sec:weight-ranks} give sums of 24 and 19
for \(p=2\) and \(p=3\), respectively. The 24 generators attain the
bound for \(p=2\). For \(p=3\) we have \(r_{k,3}=d_k\) in weights
\(40,42,46,48,50\). In each of these weights every class is a sum
of products of lower-weight classes. Eliminating the corresponding
5 generators in increasing order leaves the stated 19.
\end{proof}

\noindent\textbf{Acknowledgements.}
The author was supported by JSPS KAKENHI Grant Number 22K03259.

\medskip
\noindent\textbf{Data availability.}
Section~\ref{sec:weight-ranks} describes the rank calculations and gives the resulting ranks.

\medskip
\noindent\textbf{Conflict of interest.}
The author declares that there is no conflict of interest.

\medskip
\noindent\textbf{Use of AI.}
The author used ChatGPT for mathematical discussions and to improve the
English of the manuscript. The author takes full responsibility for the
mathematical arguments and the final manuscript.

\end{document}